\documentclass[a4paper, oneside, 10pt]{article}
\usepackage[english]{babel}
\usepackage{amsfonts}
\usepackage{amssymb}
\usepackage{graphics}
\usepackage{latexsym}
\begin{document}

\title{{\bf{\Large{On a Gibbs characterization\\ of normalized generalized Gamma processes.}}}\footnote{{\it AMS (2000) subject classification}. Primary: 60G58. Secondary: 60G09.} \footnote{{Research partially supported by MUR research grant n. 2006/134525 and Bocconi research contract n. 2007/01800.}}}
\author{\textsc {Annalisa Cerquetti}\footnote{Corresponding author. Istituto di
Metodi Quantitativi, Viale Isonzo,
25, 20133 Milano, Italy.
E-mail: {\tt annalisa.cerquetti@unibocconi.it}}\\
  \it{\small Bocconi University, Milano, Italy }}
\newtheorem{teo}{Theorem}
\date{}
\maketitle{}

\begin{abstract}
We show that a Gibbs characterization of normalized generalized Gamma processes, recently obtained in Lijoi, Pr\"unster and Walker (2007), can alternatively be derived by exploiting a characterization of exponentially tilted Poisson-Kingman models stated in Pitman (2003).
We also provide a completion of this result investigating the existence of normalized random measures inducing exchangeable Gibbs partitions of type $\alpha \in (-\infty, 0]$.\\

\end{abstract}

\section{Introduction}

The idea of constructing random probability measures by normalizing the increments of a random process dates back to Kingman's construction of the Dirichlet process by normalization of the increments of a Gamma process (Kingman, 1975). Since then it has been exploited in a variety of contexts, like e.g. ecology (Engen, 1978), population genetics (Ewens and Tavar\'e, 1995), excursion theory (Pitman and Yor, 1997) and combinatorics (Pitman, 2006). Recently the normalization approach has gained new interest in Bayesian nonparametrics, where a key problem is to define tractable alternatives to the Dirichlet process to be used as priors on the space of probability distributions (see e.g. James, 2002;  Regazzini, Lijoi and Pr\"unster, 2003). A comprehensive Bayesian analysis of a large class of random probability measures obtained by an appropriate normalization is proposed in James, Lijoi and Pr\"unster (2005) by providing explicit marginal distributions, tractable analogues of the Blackwell-MacQueen P\`olya urn prediction rules, and suitable descriptions of posterior distributions.

Normalized random measures select almost surely {\it discrete} probability measures (see e.g. James, 2003), which are well-known  to induce exchangeable random partitions of the positive integers (Kingman, 1978). Among those the subset characterized by the so-called {\it Gibbs} product form of the partition probability function, has been recently shown in Gnedin and Pitman (2006) to be a convex set, whose extreme points, for different values of the {\it type} parameter $\alpha \in (-\infty, 1)$,  have been derived by resorting to the theory of Stirling triangles. In Bayesian nonparametrics the product form of the exchangeable partition probability function (EPPF) is highly desiderable both with regard to mathematical tractability of prior-to-posterior analysis, like in connection with sequential sampling schemes for posterior's simulation, hence it is of some interest to characterize the subclass of normalized random measures which possesses this property. 

Lijoi, Pr\"unster and Walker (2007) have recently shown by means of an analytical approach based on solutions of a specific recursive equation that, in the space of (homogeneous) normalized random measures, the Gibbs product form of the EPPF for $\alpha \in (0,1)$ characterizes the class of normalized {\it generalized Gamma processes}, i.e. of random discrete probability measures whose ranked atoms follow the {\it exponentially tilted  Poisson-Kingman distribution} derived from the positive $\alpha$-stable subordinator (Pitman, 2003; see also Cerquetti, 2007).
Here we show such a characterization can also be obtained relying on known results on Poisson-Kingman models. In particular it arises by exploiting a characterization of {exponentially tilted} Poisson-Kingman models stated without proof in Pitman (2003). 
We also complete the study of normalized random measures inducing exchangeable partitions in Gibbs product form by discussing the cases $\alpha=0$ and $\alpha < 0$ which have not been treated in Lijoi, Pr\"unster and Walker (2007).

The paper is organized as follows. In Section 2 we recall basic definitions and results on normalized random measures, exchangeable random partitions and Poisson-Kingman models. In Section 3 we recall Gnedin and Pitman's (2006) characterization of Gibbs partitions in terms of extremes points and propose a detailed discussion of this fundamental result. Finally, in Section 4, after providing a proof of Pitman's characterization of exponentially tilted Poisson-Kingman models, we establish our main result. 

\section {Preliminaries and basic definitions}
We start by providing a very general and well-known construction of {\it homogeneous} normalized random measures.  
First recall that given a strictly positive random variable $T$, with  density $f_T$ and Laplace transform
\begin{eqnarray}
E(e^{-\lambda T}) =\int_{0}^{\infty} e^{-\lambda t}f_T (t)dt= e^{-\psi(\lambda)}\nonumber
\end{eqnarray}
where, according to the L\'evy-Kintchine formula, for $\lambda > 0$, $\psi(\lambda)=\int_{0}^{\infty} (1 - e^{-\lambda x})\rho(x)dx$ is the Laplace exponent, (i.e. $T$ is an infinite divisible r.v.), $f_T(\cdot)$ is uniquely identified by its unique L\'evy density $\rho(\cdot)$, which satisfies $\int_0^{\infty} \rho(x)dx=\infty$, (otherwise $P(T =0)=\exp(-\psi(\infty))>0$ contrary to the previous assumption that $T$ is strictly positive). Let $H(\cdot)$ denote a probability measure on a Polish space $(S, \mathcal {S})$, fixed and non-atomic. Now, for each $T$ and $H$, one may construct (Kingman, 1967) a {\it completely random measure} $\mu$ on $\mathcal{S}$, characterized by its Laplace functional for every positive measurable function $g$ on $S$ as
\begin{eqnarray}
\label{mualea1}
\mathbb{E}[e^{-\mu(g)}|H]=\exp\left\{-\int_{ {S}} \psi(g(s)) H(ds)\right\}=\exp\left\{-\int_{S} \int_{0}^\infty (1-e^{-g(s)x})\rho(dx) H(ds)\right\}\nonumber
\end{eqnarray}
where $\mu(g)=\int_{{S}} g(s)\mu(ds)$, so that $T=\mu({S}):=\int_{S} \mu(ds)=\int_{S} I\{s \in \mathcal {S}\}\mu(ds) $. An {\it (homogeneous) normalized random measure} (NRM) $P(\cdot)$ on $(S, \mathcal{S})$ is then obtained by normalizing $\mu$ as follows
$$
P(\cdot):=\frac{\mu (\cdot)}{\mu(S)}=\frac{\mu(\cdot)}{T}.
$$
\\
Notice that a more general construction, incorporating {\it non-homogeneous} NRMs based on non-homogeneous L\'evy measures, can be given in terms of mean intensity $\nu(dx,ds)=\rho(dx|s)H(ds)$ of a Poisson random measure $N(dx, ds)$ on $(\mathcal X \times \mathcal S)$ 
 characterized by Laplace functional
\begin{eqnarray}
\label{mualea}
\mathbb{E}[e^{-N(g)}]=\exp\left\{-\int_{S} \int_{0}^\infty (1-e^{g(s,x)}\rho(dx|s)H(ds)\right\}\nonumber,
\end{eqnarray}
(see e.g. James, 2002; James, 2005; James, Lijoi and Pr\"unster, 2005), but here we do not need to deal with such a general construction.
\\\\
As proved e.g. in James (2003) NRMs select almost surely discrete distributions, and it is well known that given a law $Q$ on the space $\mathcal{P}_1^\downarrow$ of decreasing sequences of positive numbers with sum 1, and a law $H(\cdot)$ on a Polish space $(S, \mathcal {S})$, a {\it random discrete} probability measure (RDPM) $P$ on $\mathcal{S}$ may always be defined as $P(\cdot)=\sum_{i=1}^\infty P_i \delta_{X_i}(\cdot)$, for $X_i$ iid $\sim H(\cdot)$ and $(P_i) \sim Q$. 
From Kingman's theory of exchangeable random partitions (Kingman, 1978), sampling from $P$ induces a random partition $\Pi$ of the positive integers $\mathbb{N}$ by the exchangeable equivalence relation
$i \approx j \Leftrightarrow X_i=X_j$, that is to say two positive integers $i$ and $j$ belong to the same block of $\Pi$ if and only if $X_i=X_j$, where $X_i|P$ are iid $\sim P$. It follows that, for each restriction $\Pi_n=\{A_1,\dots, A_k\}$ of $\Pi$ to $[n]=\{1,\dots, n\}$, and for each $n=1,2,\dots$, 
$$Pr(\Pi_n=\{A_1,\dots, A_k\})=p(n_1,\dots, n_k),
$$ where, for $j=1,2,\dots, k$, $n_j=|A_j|\geq 1$ and $\sum_{j=1}^k n_j=n$, for some non-negative symmetric function $p$ of finite sequences of positive integers called the {\it exchangeable partition probability function} (EPPF) determined by $\Pi$ (see Pitman, 2006, for a comprehensive account on exchangeable random partitions and related stochastic processes).\\\\ 
Pitman (2003), generalizing Kingman's (1975) construction of the Dirichlet process as a Gamma process with independent increments divided by the sum, introduces a large class of RDPMs deriving the law $Q$ 
by a random discrete distribution 
$(P_i)=(J_i/T)$, where $J_1 \geq J_2 \geq \cdots \geq 0$ are the random lenghts of the ranked points of a Poisson process with L\'evy density $\rho$ and $T=\sum_i J_i$. It is easy to see that this construction is formally equivalent to the homogeneous normalized random measure's construction given above, so that, for $(X_i)$, independent of $(P_i)$, iid $\sim H$, 
\begin{eqnarray}
\label{PKNOR}
P(\cdot)=\frac{\mu (\cdot)}{T}=\sum_{i=1}^{\infty}\frac{J_i}{T} \delta_{X_i}(\cdot)\nonumber.
\end{eqnarray}
Pitman termed the laws $Q$ of $(P_i)$ on $\mathcal{P}_1^{\downarrow}$ 
{\it Poisson-Kingman distributions with L\'evy density $\rho$}, and also enlarged the {\it basic} Poisson-Kingman family by considering the larger class  of {\it Poisson-Kingman distributions with L\'evy density $\rho$ and mixing distribution $\gamma$} given by
\begin{equation}
\label{PKMIX}
PK(\rho, \gamma):=\int_0^\infty PK(\rho|t) \gamma(dt),
\end{equation}
where $PK(\rho|t)$ is the regular conditional distribution of $(P_i)$ given $(T=t)$ constructed above, and $\gamma$ is an arbitrary probability distribution on $(0,\infty)$. Clearly if $\gamma(\cdot)=f_T(\cdot)$ then $PK(\rho,\gamma)=PK(\rho)$.\\\\  
To deal with what follows it is also worthwhile to recall that, for $(\tilde{P}_1,\tilde{P}_2,\dots)$ a {\it size-biased permutation} of the ranked atoms $(P_i)$ of $P(\cdot)$, and
$\tilde{\nu}(dp)$ the distribution of $\tilde{P}_1$ on $(0,1]$ (also termed the {\it structural distribution} of $P(\cdot)$ in Engen, 1978),  from Theorem 2.1 in Perman, Pitman and Yor (1992), (see also Lemma 1 and 2 in Pitman, 2003), the conditional distribution $PK(\rho|t)$ for $(P_i)|T=t$, 
is completely described by means of the density of the conditional {\it structural distribution} of $\tilde{P}_1|(T=t)$ given by
\begin{equation}
\label{stru}
\tilde{f}(p|t):=pt\rho(pt)\frac{f(\bar{p}t)}{f(t)}
\end{equation}
for $0<p<1$, $\bar{p}:=1-p$, $f(\cdot)$ the probability density of $T$ and $\rho(\cdot)$ the corresponding L\'evy density. This implies that the density  $\tilde{f}(p|t)$, joint with the density $\gamma(t)$, give a complete description of a mixed $PK(\rho, \gamma)$ model.\\\\ 
{\bf Example 1.} [Dirichlet process] For $\rho_\theta(x)=\theta x^{-1}e^{-x}$, $\theta>0$, the L\'evy density of the Gamma $(\theta, 1)$ density, 
with Laplace exponent $\psi(\lambda)=\theta \log(1+\lambda)$, the law $PK(\rho_\theta)$ gives the Poisson-Dirichlet ($\theta$) distribution governing the ranked atoms of the {\it normalized Gamma process}, $P(\cdot)=\mu(\cdot)/T$, for $T \sim$ Gamma$(\theta, 1)$, and $\mu(\cdot) \sim$ Gamma $(\theta H(\cdot), 1)$,  which is well-known to correspond to the Dirichlet process with parameter measure $\theta H (\cdot)$, (Ferguson, 1973; Kingman, 1975).
By (\ref{stru}) it is easy to check that structural conditional and unconditional distributions coincide and are both Beta $(1, \theta)$ with density
$\tilde{f}(p)=\theta(1-p)^{\theta -1}$,  for $0 < p< 1$. This implies $PK(\rho_{\theta})=PK(\rho_{\theta}|t)$ for every $t$, and $
PK(\rho_{\theta})=PK(\rho_{\theta}, \gamma)$, for every $\gamma$. (See Pitman, 1996)
\section{Normalized random measures inducing Gibbs partitions}
By Definition 1 in Gnedin and Pitman (2006), an exchangeable random partition $\Pi$ of the positive integers, is said to be of {\it Gibbs form} if for some nonnegative weights $W=(W_j)$ and $V=(V_{n,k})$ the EPPF of $\Pi$ can be expressed in the product form  
\begin{eqnarray}
\label{gibbs1}
p(n_1,\dots, n_k)=V_{n,k} \prod_{j=1}^k W_{n_j}\nonumber
\end{eqnarray}
for all $1 \leq k\leq n$, and all compositions $(n_1,\dots, n_k)$ of $n$. Gnedin and Pitman (2006) also show that 
to define an {\it infinite} random partition of $\mathbb{N}$, i.e. a sequence $(\Pi_n)$ consistent as $n$ varies, the weights $(W_j)$ must be of the following very special form depending on a single parameter $\alpha \in [-\infty, 1)$,
$$
W_{n_j}=(1-\alpha)_{{n_j-1}\uparrow}
$$
(with $W_j=1$ for every $j$, for $\alpha=-\infty$), and the weights $(V_{n,k})$ must be the solution to the backward recursion 
$$
V_{n,k}=(n-\alpha k)V_{n+1, k}+V_{n+1,k+1}
$$
with $V_{1,1}=1$. The solutions, for each $\alpha<1$, are then obtained identifying the extreme points of the infinite dimensional simplex of the possible weights $V$, and deriving corresponding families of extreme partitions, in terms of the laws of the corresponding ranked atoms $(P_i)$, by a combinatorial tecnique based on the theory of Stirling triangles.
The fundamental result, already stated without proof in Pitman (2003, cfr. Th. 8), is the following:\\\\
{\bf Theorem 1}. [Gnedin and Pitman, 2006; Th. 12] Each exchangeable Gibbs partition of a fixed type $\alpha \in [-\infty, 1)$, i.e. characterized by an EPPF of the form
\begin{equation}
\label{gibbsGP}
p(n_1,\dots, n_k)=V_{n,k} \prod_{j=1}^k (1-\alpha)_{n_j-1 \uparrow}
\end{equation}
is a unique probability mixture of extreme partitions of this type, which are\\\\
\begin{tabular}{ll}
a) for $\alpha \in [-\infty, 0)$ & $PD(\alpha, m|\alpha|)$ partitions with $m=0,1,\dots,\infty$,\\\\
b) for $\alpha=0$ & $PD(0,\theta)$ partitions with $\theta \in [0,\infty)$,\\\\
c) for $\alpha \in (0,1)$ &  $PK(\rho_\alpha|t)$ partitions with $t \in [0,\infty).$
\end{tabular}\\\\\\
Recall that for $0\leq\alpha<1$ and $\theta>-\alpha$, or $\alpha<0$ and $\theta=m|\alpha|$, $PD(\alpha, \theta)$ stands for the two-parameter extension of the Poisson-Dirichlet distribution (Pitman and Yor, 1997), whose general form of the EPPF is well-known to be
\begin{equation}
\label{eppfpd}
p_{\alpha, \theta}(n_1,\dots,n_k)=\frac{(\theta +\alpha)_{k-1; \alpha \uparrow}}{(1+\theta)_{n-1\uparrow}}\prod_{i=1}^{k} (1-\alpha)_{n_i-1\uparrow}.
\end{equation}
(where ${k-1; \alpha \uparrow}$ stands for the usual notation of rising factorial), and that, for $\alpha \in (0,1)$, $\rho_\alpha(\cdot)$ is the L\'evy density of the positive $\alpha$-stable distribution of index $\alpha \in (0,1)$. The EPPF induced by the conditional model $PK(\rho_\alpha|t)$ has been derived in Pitman (2003, Eq. (66)) and is given by
$$
p_{\alpha}(n_1,\dots,n_k|t)=\frac{\Gamma(1-\alpha)}{\Gamma(n -k\alpha)}\left(\frac{\alpha}{t^\alpha}\right)^{k-1}\mu_\alpha(n-1-k \alpha+\alpha|t)\prod_{i=1}^{k}(1-\alpha)_{n_1-1\uparrow}
$$
for $\mu_\alpha(q|t)=\mathbb{E}_\alpha(\tilde{P_1}^q|t)$.\\\\
Relying on Theorem 1 the problem to identify NRMs inducing EPPF of Gibbs product form (\ref{gibbsGP}), reduces to identify mixing distributions $\gamma$, respectively on $\mathbb{N}\cup\{0\}$ for case a), and on $[0,\infty)$ for case b) and c), such that the resulting partition mixture model define the ranked atoms of RDPMs which can be obtained by normalizing a completely random measures. In other words we look for mixing distributions $\gamma(\cdot)$ yielding 
corresponding strictly positive infinite divisible random variables $T$ on which to build, as from Section 2, corresponding completely random measures.  Notice that in Lijoi, Pr\"unster and Walker (2007) the discussion is confined to case c), i.e. they characterize the class of {\it normalized Generalized Gamma processes} (Pitman, 2003) as being the unique family of normalized random measures inducing Gibbs partitions of type $\alpha \in (0,1)$. Here we complete the discussion by investigating  also cases a) and b). Moreover we show that the characterization of case c) arises as a direct consequence of a characterization of {\it exponentially tilted} Poisson-Kingman models stated in Pitman (2003) for which we also provide a proof.\\\\
\noindent{\bf Case a)} For $\alpha \in [-\infty, 0)$ a law $Q$ on $\mathcal{P}_1^\downarrow$ induces an EPPF of Gibbs form (\ref{gibbsGP}) if  
\begin{eqnarray}
\label{fisher}
Q_{\alpha, \gamma}(\cdot)=\sum_{m=0}^\infty PD(\alpha, m|\alpha|)\gamma(m)\nonumber
\end{eqnarray}
for $\gamma(\cdot)$ a probability distribution on the space of non negative integers. Recall that, for each $m$, $PD(\alpha, m|\alpha|)$ are Fisher's models for species sampling (Fisher et al. 1943, see also Pitman, 1996) and correspond to symmetric Dirichlet random vectors of dimension $m$ and parameter $|\alpha|$. This implies that a random $P$ whose ranked atoms are $PD(\alpha, m|\alpha|)$ distributed can be constructed as
\begin{eqnarray}
\label{findir}
P(\cdot)=\sum_{i=0}^m \frac{G_i}{G} \delta_{X_i}(\cdot)\nonumber
\end{eqnarray}
where $G_i$ are iid Gamma $(-\alpha, 1)$ r.v.'s, $G=\sum_{i=1}^m G_i$ is Gamma $(-\alpha m, 1)$, and $(X_i)$ are iid $\sim H$ independent of $(G_i)$. In Bayesian nonparametric literature $P(\cdot)$ is also called a {\it finite dimensional Dirichlet prior} (see e.g. Ishwaran and Zarepour, 2002). Mixing with respect to $m$ corresponds to randomize the number $m$ of terms in $G$, so that $G=\sum_{i=0}^M G_i$, for $M \sim \gamma(\cdot)$, has {\it compound} $\gamma-${\it Gamma distributions}, with a positive mass at zero. To our purpose we want $G$ to be infinite divisible (ID) with L\'evy measure having  infinite total mass.  From e.g. Theorem 3.2 in Steutel and Van Harn (2004) a compound distribution with positive mass at zero is ID if and only if it is a {\it compound Poisson}, hence we need to restrict our search to $\gamma$ having Poisson distribution. Recall that compound Poisson-Gamma distributions (Aalen, 1992) belong to the family of {\it generalized Gamma} distributions (see e.g. Brix 1999), 
a class of infinite divisible distributions defined for $\alpha <1$ and characterized by Laplace exponent of the form
$$
\psi_{GG}(\lambda)= -\delta[\zeta + (\zeta^{\frac 1\alpha} +2\lambda)^{\alpha}]
$$
for $\lambda \geq 0$, $\delta >0$ and $\zeta \geq 0$. Neverthless, for $\alpha \in (-\infty, 0)$, $\psi_{GG}(\infty)=\int_0^\infty \rho_{GG}(dx)< \infty$, i.e. $P(T=0)>0$, (cfr. also Pitman, 2003, Sec. 5.2), therefore no random probability measure whose ranked atoms induce Gibbs partitions of type $\alpha <0$, can be obtained by normalizing a completely random measure.\\\\  
{\bf Case b)} For $\alpha=0$ a law $Q_{\gamma, \theta}$ on $\mathcal{P}_1^\downarrow$ induces an EPPF of Gibbs form (\ref{gibbsGP}) if 
$$
Q_{\theta, \gamma}(\cdot)=\int_0^\infty PD(0,\theta)\gamma(d\theta) 
$$
where $PD(0, \theta)=PK(\rho_\theta)$ is the Poisson-Dirichlet $(\theta)$ distribution of the ranked atoms of the normalized Gamma process $P(\cdot)=\frac{\mu(\cdot)}{T}$, for $T\sim G(\theta, 1)$ and $\mu(\cdot) \sim G(\theta H(\cdot),1)$ as recalled in Example 1. From (\ref{eppfpd}) the EPPF of $PD(0, \theta)$ arises from letting $\alpha=0$ and it is in Gibbs product form. About mixing over $\theta$ notice that it corresponds to mixing over the total mass of the parameter measure  $\theta H(\cdot)$ of the Dirichlet process, 
therefore $Q_{\theta, \gamma}$ provides the distribution of the ranked atoms of {\it mixtures of Dirichlet processes} as introduced in Antoniak (1974). 

It is well-known that this kind of mixing 
induces correlation (see e.g. Sibisi and Skilling, 1997, Sec. 11) i.e. the resulting random process has no longer independent increments, hence it no longer corresponds to a completely random measure. This implies that, apart from the normalized Gamma process, whose ranked atoms are $PD(0, \theta)$ distributed, no other normalized random measure can induce exchangeable partitions in Gibbs product form of type $\alpha=0$.\\\\
\noindent{\bf Case c)} For $\alpha \in (0,1)$ a law $Q_{\alpha, \gamma}$ on $\mathcal{P}_1^\downarrow$ induces an EPPF of Gibbs form (\ref{gibbsGP}) if  
\begin{equation}
\label{dodici}
Q_{\alpha, \gamma}(\cdot)=\int_0^\infty PK(\rho_{\alpha,\delta}|t)\gamma(t)dt
\end{equation}
for $\rho_{\alpha, \delta}$, the L\'evy density of the positive $(\alpha,\delta)$--stable density $f_{\alpha, \delta}$
\begin{eqnarray}\rho_{\alpha,\delta}(x)=\delta2^{\alpha}\frac{\alpha}{\Gamma(1-\alpha)}x^{-1-\alpha}\nonumber,
\end{eqnarray}
and $\gamma(\cdot)$ an arbitrary probability density on $(0, \infty)$. Hence the laws $Q_{\alpha, \gamma}$ correspond to mixed Poisson-Kingman models $PK(\rho_\alpha, \gamma)$, as introduced in Section 2, derived from the positive stable law of index $\alpha \in (0,1)$ and scale parameter $\delta$. This is the case treated in Lijoi, Pr\"unster and Walker (2007, cfr. Prop. 2). 
For a general $\alpha \in (0,1)$ the conditional structural distribution of $_\alpha\tilde{P}_1|T=t$ which describes $PK(\rho_{\alpha, \delta}|t)$, has density
\begin{eqnarray}
\label{cond}\tilde{f}_{\alpha, \delta}(p|t)=\left(\frac{2}{pt}\right)^{\alpha}\frac{\delta\alpha}{\Gamma(1-\alpha)}\frac{f_{\alpha,\delta}(\bar{p}t)}{f_{\alpha,\delta}(t)}
\end{eqnarray}
but has no closed form expression since explicit expressions of $f_{\alpha, \delta}$ are known only in the form of series representation. For $\alpha=\frac 12$ the stable density corresponds to the L\'evy density, 
$f_{\frac12, \delta}(t)=\frac{\delta}{\sqrt{2\pi}}t^{-\frac 32} \exp\{-\frac {\delta^2}{2t}\}$
and by substitution in (\ref{cond})
\begin{eqnarray}
\label{csl} \tilde{f}_{\frac 12,\delta}(p|t)=\frac{\delta}{\sqrt{2\pi}\sqrt{pt}}(1-p)^{-\frac 32}\exp\left\{-\frac{1}{2}\frac{p\delta^2}{(1-p)t}\right\}.
\end{eqnarray}
For $\gamma(t)=f_{\alpha, \delta}(t)$, $PK(\rho_{\alpha,\delta}, f_{\alpha, \delta})=PK(\rho_{\alpha,\delta})$ gives the law of the ranked atoms of a {\it normalized $\alpha$--stable} process, with Laplace exponent $\psi_{\alpha,\delta}(\lambda)=\delta(2\lambda)^{\alpha}$.\\\\
Now notice that, in general, 
apart from the trivial case $\gamma(t)=f_T(t)$, {\it mixed} $PK(\rho, \gamma)$ models describe the law of the ranked atoms of a random discrete probability measure which can be obtained by normalization of a 
completely random measure $P(\cdot)$, if only if admit an equivalent formulation as {\it basic} Poisson-Kingman models.
The following proposition can be stated that doesn't need to be proved:\\\\
{\bf Proposition 1.} Let $\rho(\cdot)$ be a L\'evy density on $(0,\infty)$, corresponding to a strictly positive infinitely divisible r.v. $T$, let 
$\gamma \neq f_T$ be a probability density on $(0, \infty)$  and $H$ a fixed distribution on $(S, \mathcal{S})$ then, a {\it mixed} $PK(\rho, \gamma)$ model on $\mathcal{P}_1^{\downarrow}$ gives the law of the ranked atoms of a random discrete probability measure $P(\cdot)=\sum_i P_i \delta_{X_i}(\cdot)$ where $(X_i)$, independent of $(P_i)$, are iid $\sim H$,  which admits a construction by normalization of a completely random measure if and only if there exists a Levy density $\rho^*(\cdot)$ with $\int_0^\infty \rho^*(x)dx=\infty$, such that
\begin{equation}
\label{equi}
PK(\rho^*)\equiv PK(\rho, \gamma).\\\\
\end{equation}
Relying on Proposition 1, in case c) our problem reduces to identify the family of mixing distributions $\gamma(\cdot)$ which satisfy equation (\ref{equi}) for $\rho=\rho_\alpha$. To this end, in the following section, we resort to a result stated in Pitman (2003) for which we also provide a proof.\\\\
\section {Pitman's characterization of exponentially tilted PK models.}
In Section 4.2 Pitman (2003) focuses on {\it exponential tilting} as one of the basic operations on L\'evy densities which lead to a tractable class of mixed Poisson-Kingman partitions models.  The idea of tilting density functions is very old, in L\'evy processes setting the equivalent transformation it is also known as {\it Esscher transform} (see e.g. Sato, 1999). Here we recall the basic definition:\\\\
{\bf Definition 1.} [Exponential tilting] Given a probability density $f$ on $(0,\infty)$, with Laplace exponent $\psi(\lambda)$, for every $\lambda >0$, the corresponding {\it exponentially tilted density} $f_\lambda$ is obtained as
\begin{eqnarray}
f_\lambda(t)=\exp\{-\lambda t + \psi(\lambda) +k(t)\}=\frac{1}{L(\lambda)}e^{-\lambda t} f(t)\nonumber
\end{eqnarray}
where $f(t)=\exp\{k(t)\}$ and $L(\lambda)=E(e^{-\lambda T})=\exp\{-\psi(\lambda)\}$. The corresponding 
Laplace transform is given by
\begin{equation}
\label{laptilt}
\exp\{-\psi_\lambda(b)\}=\exp\{-\psi(b+\lambda)+\psi(\lambda)\},
\end{equation}
for $b>0$, if additionally $f$ is {\it infinitely divisible}, then
\begin{eqnarray}
\psi_\lambda(b)=\int_0^\infty (1-e^{-bs})e^{-\lambda s} \rho(ds),\nonumber
\end{eqnarray}
hence tilting $f$ yields corresponding tilted L\'evy measures, $\rho_\lambda(\cdot)=e^{-\lambda t}\rho(\cdot)$.\\\\
What Pitman (2003) basically states is that {\it exponential tilting} is the unique operation on basic PK models that produces mixed PK models satisfying condition of Proposition 1. More exactly a {\it basic} model equivalent to a {\it mixed}  model exists if and only if the mixing density belongs to the family of the corresponding tilted densities, in which case the basic model is driven by the tilted version of the L\'evy measure of the mixed model. (Notice that the ``only if'' part of the statement is partially hidden in the last line of Section 4.2 in Pitman's paper). 

Here we recall Pitman's result, both elucidating the ``if'' part and providing a proof for the ``only if'' part.\\\\
{\bf Proposition 2.} [Pitman, 2003; Sec. 4.2] Given two regular L\'evy densities $\rho_2$ and $\rho_1$, then $$PK(\rho_2)=PK(\rho_1, \gamma)$$ for some $\gamma(\cdot)$ if and only if $\rho_2(t)=\rho_1(t) \exp\{-\lambda t\}$ and $\gamma(t)=f_1(t)\exp\{\psi_1(\lambda) -\lambda t\}$, for every $t$, for some $\lambda >0$, where $f_1(\cdot)$ is the probability density corresponding to $\rho_1$ and $\psi_1(\lambda)$ its Laplace exponent.\\\\ 
{\it Proof}
{(\it Sufficiency)} Let $\rho(\cdot)$ be a L\'evy density corresponding to a probability density $f(\cdot)$, and $\rho_\lambda(\cdot)=e^{-\lambda t}\rho(\cdot)$ the tilted L\'evy density of the corresponding tilted probability density $f_\lambda(\cdot)=f_1(t)\exp\{\psi_1(\lambda) -\lambda t\}$. 
It is easy to show that the conditional structural distributions for the models $PK(\rho|t)$ and $PK(\rho_\lambda|t)$ coincide, in fact 
$$
\tilde{f}_\lambda(p|t)=pt\rho_\lambda(pt)\frac{f_\lambda(\bar{p}t)}{f_\lambda(t)}= pte^{-\lambda pt}\rho(pt)\frac{e^{\psi(\lambda)-\lambda \bar{p}t}f(\bar{p}t)}{e^{\psi(\lambda)-\lambda t}f(t)}=pt\rho(pt)\frac{f(\bar{p}t)}{f(t)}=\tilde{f}(p|t).
$$
It follows that $PK(\rho_\lambda|t)=PK(\rho|t)$
for every $t$, hence 
\begin{eqnarray}
\label{picappa}
PK(\rho_\lambda, \gamma)=PK(\rho, \gamma)\nonumber
\end{eqnarray}
for every mixing density $\gamma$. For $\gamma_\lambda=f_\lambda=f_1(t)\exp\{\psi_1(\lambda) -\lambda t\}$, this implies $PK(\rho_\lambda, \gamma_\lambda)=PK(\rho, \gamma_\lambda)$, 
hence
$$
PK(\rho_\lambda)=PK(\rho_\lambda, \gamma_\lambda)=PK(\rho, \gamma_\lambda).
$$
\hspace{15cm}$\square$\\
({\it Necessity}) 
We want to show that if $\rho_1$ is the L\'evy density of some infinitely divisible density $f_1(\cdot)$ and there exists a L\'evy density $\rho_2(x)\neq \rho_1(x)$ such that, for some density $\gamma$
\begin{equation}
\label{equiv}
PK(\rho_2)=PK(\rho_1, \gamma),
\end{equation}
then $\rho_2(\cdot)=e^{-\lambda t}\rho_1(\cdot)$ and $\gamma(t)=e^{\psi_1(\lambda)-\lambda t}f_1(t)$, for $\lambda >0$ and $\psi_1(\cdot)$ the Laplace exponent of $f_1(t)$. 
Now, call $f_2$ the probability density corresponding to the L\'evy density $\rho_2$ of some infinitely divisible r.v. $T_2$. Since $PK(\rho_2)=PK(\rho_2, f_2)$, then (\ref{equiv}) equates 
\begin{equation}
\label{equivdue}
PK(\rho_2, f_2)=PK(\rho_1, \gamma).
\end{equation}
Condition (\ref{equivdue}) implies that all the following three conditions must hold:\\\\
a) there exists some function $\phi(\cdot)$ such that $\gamma(t)=\phi(t)f_1(t)=f_2(t)$ for every $t$,\\\\
b) there exists a function $g(\cdot)$ on $(0, \infty)$  such that, for $b>0$ 
$$
E(e^{-\lambda T_2})=\int_0^{\infty} e^{-b t} \phi(t)f_1(t)dt=\exp\{-\psi_2(b)\}=\exp{-\int_0^\infty (1-e^{-b x})g(x)\rho_1 (dx)}\\\\
$$
is the Laplace transform of $T_2$,\\\\
c) $PK(\rho_2|t)=PK(\rho_1|t)$, for every $t$.\\\\
Condition b) can be rewritten as 
\begin{eqnarray}
\label{intuno}
\int_0^{\infty} e^{\psi_2(b)}e^{-b t} \phi(t)f_1(t)dt=1\nonumber
\end{eqnarray}
for which a family of possible solutions for every $b>0$ is given by:
\begin{equation}
\label{formaphi}
\phi(t, b)=e^{-\psi_2(b)+b t}.
\end{equation} 
Conditions b) and c) imply
$$
\rho_1(pt)\frac{f_1(\bar{p}t)}{f_1(t)}=\rho_1(pt)g(pt)\frac{f_1(\bar{p}t)\phi(\bar{p}t)}{f_1(t)\phi(t)}
$$
for all $t$ and for all $0<p<1$, which yields 
\begin{equation}
\label{phiexp}
\phi(t)=\phi(\bar{p}t)g(pt).
\end{equation} 
For the uniqueness of the Laplace transform, the density $f_2(t)$ corresponding to $\psi_2$ must be unique, therefore combining conditions (\ref{formaphi}) and (\ref{phiexp}) one obtains
$$
e^{-\psi_2(b)+b t}=e^{-\psi_2(b)+b \bar{p}t}g({p}t)
$$
which yields $g(pt)=e^{b p t}$ for $b>0$. It follows that the  transformed L\'evy density must be of the form $\rho_2=e^{b t}\rho_1$. By Example 33.15 in Sato (1999), for $b<0$, $g(x)=e^{bx}$ is the {\it exponential transformation}, hence by the change of variable $\lambda=-b$, it follows that $\rho_2=e^{-\lambda t}\rho_1(\cdot)$  and the corresponding  tilted probability density 
$$
f_2(t)=e^{\psi_1(\lambda)-\lambda t}f_1(t)
$$
are the unique solutions to (\ref{equiv}).\hspace{10cm}$\square$\\\\
{\bf Remark  1.} [Dirichlet process] By Example 1, $PK(\rho_\theta)=PK(\rho_\theta, \gamma)$
for every $\gamma$, hence $PK(\rho_\theta)=PK(\rho_\theta, \gamma_{\theta, \lambda})$
for $\gamma_{\theta, \lambda}$ the tilted Gamma density. Moreover, by Proposition 2, $PK(\rho_{\theta,\lambda})=PK(\rho_\theta, \gamma_{\theta, \lambda})$, for $\rho_{\theta,\lambda}$ the corresponding L\'evy measure, hence $PK(\rho_\theta)=PK(\rho_{\theta,\lambda})$ for every $\lambda$. 

It is easy  in fact to verify that the operation of tilting a Gamma $(\theta, 1)$ density 
$$
f_\lambda(t)=\frac{1}{\Gamma(\theta)}x^{\theta-1}e^{-x} exp\{-\theta log(1+\lambda)-\lambda x\}=\frac{(1+\lambda)^{\theta}}{\Gamma(\theta)}x^{\theta-1}e^{-x(1+\lambda)}
$$
yields a family of Gamma$(\theta, (1+\lambda))$ laws.\\\\
We are now in a position to complete the discussion of case c) of Theorem 1.
By exploiting Proposition 2, the solution to case c) arises immediately. The law $Q_{\alpha, \gamma}(\cdot)$ in (\ref{dodici}) admits a representation as a basic PK models, if and only if  $\gamma$ belongs to the family of the {\it exponentially tilted positive stable} densities, (first introduced in Hougaard, 1986, also called {\it tempered stable} distributions in Barndorff-Nielsen and Shepard, 2001) given by, 
$$f_{\alpha, \delta}^{\lambda}(t )=f_{\alpha,\delta}(t)\exp\{\delta \lambda^{\alpha}-\lambda t\}.$$
The basic model is driven by the corresponding tilted L\'evy density
$$
\rho_{\alpha, \delta}^{\lambda}(s)=\frac{\delta2^{\alpha}\alpha}{\Gamma(1-\alpha)}s^{-1-\alpha}e^{-\lambda s}.
$$  
By the change of variable $\lambda=\frac{\zeta^{\frac 1\alpha}}{2}$, and applying (\ref{laptilt}), the Laplace exponent of $f_{\alpha, \delta}^{\lambda}$results
\begin{equation}
\label{lapl}
\psi^{\zeta}_{\alpha,\delta}(b)=-\delta\zeta+\delta(\zeta^{\frac 1\alpha}+2b)^{\alpha}
\end{equation}
which is well-known to identify to the family of infinitely divisible {\it generalized Gamma distributions} (Brix, 1999,) as already recalled in the discussion of case a). Notice that, contrary to case $\alpha < 0$, for $\alpha \in (0,1)$, $\psi_{GG}(\infty)=\infty$, hence $P(T>0)=1$.\\\\
{\bf Example 2.} [Inverse-Gaussian process]  
As previously recalled for $\alpha=\frac 12$ the positive {\it stable} density has explicit expression
$$
f_{\frac12, \delta}(t)=\frac{\delta}{\sqrt{2\pi}}t^{-\frac 32} e^{-\frac {\delta^2}{2t}}\nonumber, 
$$
with corresponding Laplace exponent $\psi(\lambda)=\delta\sqrt{2\lambda}$ and L\'evy density $\rho_{\frac 12, \delta}=\frac{\delta}{\sqrt{2\pi}}x^{-\frac 32}\nonumber$. By exponential tilting with $\lambda=\frac {\zeta^2}{2}$, the {\it exponentially tilted} $1/2-${\it stable} density results 
\begin{equation}
\label{igdens}
f_{\delta, \zeta}(t)= 
\frac{\delta}{\sqrt{2\pi}} e^{\delta \zeta}{t^{-\frac 32}} 
\exp\left\{-\frac 12 \left(\delta^2t^{-1}+\zeta^2 t\right)\right\},
\end{equation}
for $\delta \in (0,\infty)$ and $\zeta \in [0,\infty)$, which is well-known to be the density of an {\it inverse Gaussian $(\delta, \zeta)$ law} (see e.g. Seshadri, 1993). By (\ref{lapl}) corresponding Laplace exponent results 
$
\psi_{\frac 12, \delta}^\zeta(b)=-\delta\zeta +\delta(\zeta^2 +2b)^{\frac 12}
$
and the corresponding tilted L\'evy density is given by:
$$
\rho_{\frac 12, \delta}^\zeta(x)=\frac {\delta}{\sqrt{2\pi}}x^{-\frac 32}e^{-\frac x2 \zeta^2}.
$$
By Proposition 2 
\begin{eqnarray}
\label{PKNIG}
PK(\rho_{\frac 12, \delta}^\zeta)=PK(\rho_{\frac 12, \delta},f_{\delta, \zeta})\nonumber,
\end{eqnarray}
and it is an easy task to verify that the conditional structural distribution derived from the inverse Gaussian density equates (\ref{csl}). Applying (\ref{stru}) for $f_{\alpha, \zeta}$ as in (\ref{igdens})
$$
\tilde{f}_{\delta, \zeta}(p|t)= \frac {pt\delta}{\sqrt{2\pi}}(pt)^{-\frac 32}\exp\left\{-\frac {pt}{2} \zeta^2\right\}\frac{{(pt)^{-\frac 32}} 
\exp\left\{-\frac 12 \left(\delta^2((1-p)t)^{-1}+\zeta^2 (1-p)t\right)\right\}}{{t^{-\frac 32}} 
\exp\left\{-\frac 12 \left(\delta^2t^{-1}+\zeta^2 t\right)\right\},}
$$ 
and some elementary calculations show it simplifies to 
$$
\tilde{f}_{IG}(p|t)=\frac{\delta}{\sqrt{2\pi}\sqrt{pt}}(1-p)^{-\frac 32}\exp\left\{-\frac{1}{2}\frac{p\delta^2}{(1-p)t}\right\}.
$$
Notice that the corresponding normalized random measure, termed {\it normalized Inverse Gaussian process} 
has been derived by mimicking Ferguson's (1973) construction of the Dirichlet process, and studied in relation to hierarchical Bayesian nonparametric mixture modeling in Lijoi, Mena and Pr\"unster (2005).\\\\
This complete the discussion of case a), b), and c) of Gnedin and Pitman's characterization of Gibbs partitions. We have shown that for $\alpha \in (-\infty, 0)$ nor the extreme partitions, nor the family of mixture partition models, contain elements that can be derived by means of normalization of completely random measures. For $\alpha=0$, only the family of extreme partitions $\{PD(0,\theta), \theta >0\}$, correspond to the ranked atoms of a NRM, namely the {\it normalized Gamma} process, which belongs to the family of Generalized Gamma processes for $\alpha=0$. Finally for $\alpha \in (0,1)$ only mixing exponentially tilted  stable distributions of index $\alpha \in (0,1)$ produce partitions models for random discrete probability measures that can be derived by normalization of completely random measures, namely {\it generalized Gamma} processes of index $\alpha \in (0,1)$. \\\\
We summarize our conclusions in the following proposition, slightly more general than that given in Lijoi, Pr\"unster and Walker (2007, Prop. 2).\\\\ 
{\bf Proposition 3.} The unique family of random discrete probability measures which admits a construction through normalization of completely random measures, and induces exchangeable partition probability functions in Gibbs product form of type $\alpha <1$ is the class of normalized generalized Gamma processes of index $\alpha \in [0,1)$.\\\\ 
In Cerquetti (2007) an explicit form of the EPPF of 
exponentially tilted Poisson-Kingman models derived from the positive stable subordinator for $\alpha \in (0,1)$, that here we recall for completeness, is derived from the general form for mixed PK models given in Pitman (2003), (see also Lijoi, Mena and Pr\"unster, 2007),
\begin{eqnarray}
\label{EPPFPKA}
p_{GG}(n_1,\dots,n_k)=\frac{e^{\delta\gamma} \delta^k \alpha^k 2^n }{\Gamma(n)} \prod_{j=1}^k (1-\alpha)_{n_j -1 \uparrow}\int_{0}^{\infty} \lambda^{n-1} \frac{e^{-\delta (\gamma^{\frac 1\alpha}+2\lambda)^\alpha}}{(\gamma^{\frac 1\alpha}+2\lambda)^{n-k\alpha}}d\lambda.\nonumber
\end{eqnarray}
Notice the Gibbs product form of type $\alpha$ of (\ref{EPPFPKA}), where
\begin{eqnarray}
&&V_{n,k} =\frac{e^{\delta\gamma} \delta^k \alpha^k 2^n }{\Gamma(n)} \int_{0}^{\infty} \lambda^{n-1} \frac{e^{-\delta (\gamma^{\frac 1\alpha}+2\lambda)^\alpha}}{(\gamma^{\frac 1\alpha}+2\lambda)^{n-k\alpha}}d\lambda\nonumber.
\end{eqnarray}

\section*{References}
\newcommand{\bib}{\item \hskip-1.0cm}
\begin{list}{\ }{\setlength\leftmargin{1.0cm}}

\bib \textsc{Aalen, O. O.} (1992) Modelling heterogeneity in survival analysis by the compound Poisson distribution. {\it Ann. Appl. Probab.} 2,  951-972.

\bib \textsc{Antoniak, C. E.} (1974) Mixtures of Dirichlet processes with applications to Bayesian nonparametric problems. {\it Ann. Statist.} 2, 1152-1174. 

\bib \textsc{Barndorff-Nielsen, O. E. and Shepard, N.} (2001)  Normal modified stable processes.  {\it Th. Probab. Math. Statist.}, 65, 1-19.

\bib \textsc{Brix, A.} (1999) Generalized Gamma measures and shot-noise Cox processes. {\it Adv. Appl. Probab.}, 31, 929--953.

\bib \textsc{Cerquetti, A.} (2007) A note on Bayesian nonparametric priors derived from exponentially tilted Poisson-Kingman models. {\it Stat. \& Prob. Lett}, (To appear).

\bib \textsc{Engen, S.} (1978) {\it Stochastic abundance models}. Chapman \& Hall, London.

\bib \textsc{Ewens, W. and Tavar\'e S.} (1995) The Ewens sampling formula. In Multivariate discrete distributions (Johnson, N.S., Kotz, S. and Balakrishnan, N. eds.). Wiley, NY.

\bib \textsc{Ferguson, T. S.} (1973)  A Bayesian analysis of some nonparametric problems. {\it Ann. Statist.}, 1, 209--230.

\bib \textsc{Fisher, R.A., Corbet, A.S. and Williams, C.B.} (1943) The relation beteween the number of species and the number of individuals in a random sample of an animal population. {\it J. Animal. Ecol.}, 12, 42--58.

\bib \textsc{Gnedin, A. and Pitman, J. } (2006) {Exchangeable Gibbs partitions  and Stirling triangles.} {\it Journal of Mathematical Sciences}, 138, 3, 5674--5685. 



\bib \textsc{Hougaard, P} (1986) Survival models for hetereneous populations derived from stable distributions. {\it Biometrika}, 73, 387--396.

\bib \textsc{Ishwaran, H. \& Zarepour, M.} (2002) Exact and approximate sum-represen\-tations for the Dirichlet process. {\it Can. J. Statist.} 30, 269-283. 

\bib \textsc{James, L. F.} (2002). Poisson process partition calculus with applications to exchangeable models and Bayesian Nonparametrics. {\it arXiv:math.ST/0205093}.

\bib \textsc{James, L. F.} (2003) A simple proof of the almost sure discreteness of a class of random measures. {\it Statist. \& Probab. Lett.}, 65, 363-368.

\bib \textsc{James, L.F., Lijoi, A. and Pr\"unster I.} (2005) Bayesian inference via classes of normalized random measures. {\it arXiv:math.ST/0503394}.

\bib \textsc{Kingman, J.F.C} (1967) Completely random measures. {\it Pacific J. Math.}, 21, 59-78

\bib \textsc{Kingman, J.F.C.} (1975) Random discrete distributions. {\it J. Roy. Statist. Soc. B}, 37, 1--22. 

\bib \textsc{Kingman, J.F.C} (1978) The representation of partition structure.  {\it J. London Math. Soc.} 2, 374--380.

\bib \textsc{Lijoi, A., Mena, R. and Pr\"unster, I.} (2005) Hierarchical mixture modeling with normalized Inverse-Gaussian priors. {\it JASA}, vol. 100, 1278--1291.

\bib \textsc{Lijoi, A., Mena, R. and Pr\"unster, I.} (2007) Controlling the reinforcement in Bayesian nonparametric mixture models. {\it J. Roy. Statist. Soc. B}, (To appear).

\bib \textsc{Lijoi, A., Pr\"unster, I. and Walker, S.G.} (2007)  Investigating nonparametric priors with Gibbs structure. {\it Statistica Sinica}, (To appear).

\bib \textsc{Perman, M., Pitman, J, \& Yor, M.} (1992) Size-biased sampling of Poisson point processes and excursions. {\it Probab. Th. Rel. Fields}, 92, 21--39.

\bib \textsc{Pitman, J.} (1996) Some developments of the Blackwell-MacQueen urn scheme. In T.S. Ferguson, Shapley L.S., and MacQueen J.B., editors, {\it Statistics, Probability and Game Theory}, volume 30 of {\it IMS Lecture Notes-Monograph Series}, pages 245--267. Institute of Mathematical Statistics, Hayward, CA.

\bib \textsc{Pitman, J.} (2003) {Poisson-Kingman partitions}. In D.R. Goldstein, editor, {\it Science and Statistics: A Festschrift for Terry Speed}, volume 40 of Lecture Notes-Monograph Series, pages 1--34. Institute of Mathematical Statistics, Hayward, California.

\bib \textsc{Pitman, J.} (2006) {\it Combinatorial Stochastic Processes}. Ecole d'Et\'e de Probabilit\'e de Saint-Flour XXXII - 2002. Lecture Notes in Mathematics N. 1875, Springer.

\bib \textsc{Pitman, J. and Yor, M.} (1997) The two-parameter Poisson-Dirichlet distribution derived from a stable subordinator. {\it Ann. Probab.}, 25:855--900.

\bib \textsc{Regazzini, E., Lijoi, A. and Pr\"unster, I.} (2003) Distributional results for means of random measures with independent increments. {\it Ann. Statist.}, 31, 560--585.

\bib \textsc{Sato, K.} (1999) {\it L\'evy processes and infinitely divisible distributions}. Cambridge University Press.

\bib \textsc{Seshadri, V. (1993)} {\it The inverse Gaussian distribution}. Oxford University Press, New York. 

\bib \textsc{Sibisi, S. \& Skilling, J.} (1997) Prior distributions on measure space. {\it J.R. Statist. Soc.} B, 59, 1, 217--235.

\bib \textsc{Steutel F, W, \& Van Harn, K.} (2004) {\it Infinite divisibility of probability distributions on the real line}. {Monographs and Textbooks in Pure and Applied Mathematics}, 259. Marcel Dekker, New York.

\end{list}

\end{document}